\documentclass{article}
\usepackage{ijcai17}

\usepackage{times}

\usepackage{amssymb}
\usepackage{mathtools}
\usepackage{amsmath}
\usepackage{gensymb}
\usepackage{algorithm}
\usepackage{algorithmic}
\usepackage{graphicx}
\usepackage{graphicx}
\usepackage[utf8]{inputenc}
\usepackage{listings}
\usepackage{color}
\usepackage{upquote}
\usepackage{titlesec}
\usepackage{lipsum}
\usepackage{float}
\usepackage{mathtools}
\usepackage{caption}
\usepackage{subcaption}
\usepackage{amsmath}
\usepackage{relsize}
\usepackage{graphicx}
\usepackage[utf8]{inputenc}
\usepackage{listings}
\usepackage{color}
\usepackage{upquote}
\usepackage{titlesec}
\usepackage{lipsum}
\usepackage{float}
\usepackage{mathtools}
\usepackage{caption}
\usepackage{subcaption}

\definecolor{dkgreen}{rgb}{0,0.6,0}
\definecolor{gray}{rgb}{0.5,0.5,0.5}
\definecolor{mauve}{rgb}{0.58,0,0.82}

\title{A comparison between the Split Step Fourier and Finite-Difference method in analysing the soliton collision of a type of Nonlinear Schrödinger equation found in the context of optical pulses}
\author{Luke Taylor\\ 
University of Cape Town\\
Department of Mathematics and Applied Mathematics\\
Cape Town, South Africa\\
tylchr011@uct.ac.za}

\begin{document}

\maketitle

\begin{abstract}
In this report a type of Schrödinger Equation which is found in the context of optical pulses is analysed using the \textit{Split Step} and \textit{Finite Difference} method. The investigation shows interesting dynamics regarding certain values for parameter $S$ as well as a comparison between the two numeric schemes demonstrating the \textit{Split Step} to be superior for this problem.
\end{abstract}

\section{Introduction}

Optical pulses in media with saturation nonlinearity properties are modelled by the nonlinear Schrödinger \textit{(NLS)} Equation \cite{nls} with the following nonlinearity \cite{nlsoptics}. \\
\begin{equation}
\label{eq:1}
\centering \displaystyle i\dfrac{\partial \psi}{\partial t} + \dfrac{1}{2} \dfrac{\partial^2 \psi}{\partial x^2} + \dfrac{|\psi|^2\psi}{1+S|\psi|^2} = 0
\end{equation}\\
This is a Partial Differential Equation \textit{(PDE)} as it describes a relation of $\psi$ in regards to change over time and space. A solution to this equation is a so called soliton of the form.
\begin{equation}
\label{eq:2}
\centering \displaystyle \psi(x, t) = \dfrac{2 \sqrt{2}e^{\sqrt{2}x}}{1 + (\frac{3}{2}-2S)e^{2\sqrt{2}x}}e^{it+iv}
\end{equation}\\
The aim of this report is to study the collision of two such solitons \cite{soliton} with varying values for $S$ and $v$ by numerically advancing from an initial configuration (the sum of two solitons) via the \textit{Split-Step} \cite{splitstep} and \textit{Finite Difference} \cite{finitedifference} methods.

\section{Split Step Method}
The \textit{Split Step} method is a pseudo-spectral numerical method used to solve nonlinear \textit{PDEs} such as \textit{PDE} \ref{eq:1}. The method works by splitting the equation into a \textit{nonlinear} and \textit{linear} part.
\begin{equation}
\label{eq:3}
\centering \displaystyle i\dfrac{\partial \psi}{\partial t}  + \dfrac{|\psi|^2\psi}{1+S|\psi|^2} = 0
\end{equation}\\
\begin{equation}
\label{eq:4}
\centering \displaystyle i\dfrac{\partial \psi}{\partial t} + \dfrac{1}{2} \dfrac{\partial^2 \psi}{\partial x^2} = 0
\end{equation}\\

Both these equations are treated separately. The soliton is advanced in time by taking a small time step $\tau$ for both solutions. For the linear solution however $\psi$ needs to be Fourier transformed with the solution being advanced in Fourier space before inverse Fourier transforming back to the time domain.

\subsection{Solving the nonlinear part}

In order to solve equation \ref{eq:3} we multiply both sides by $i$ which gives us.
\begin{equation}
\label{eq:5}
\centering \displaystyle \dfrac{\partial \psi}{\partial t} = i\dfrac{|\psi|^2}{1+S|\psi|^2}\psi
\end{equation}\\
Next we can notice that $|\psi|^2$ is a scalar and hence equation \ref{eq:5} is a first order differential equation which has the following analytical solution.
\begin{equation}
\label{eq:5}
\centering \displaystyle \psi(x,t)=\psi(x,t_0)e^{\frac{|\psi|^2}{1+S|\psi|^2}t}
\end{equation}\\

\subsection{Solving the linear part}

The key insight is to write $\psi(x,t)$ as a Fourier Series.

\begin{equation}
\label{eq:6}
\centering \displaystyle \psi(x,t)= \sum\limits_{i=-\infty}^{\infty} \hat{\psi}_n(t)e^{\frac{2\pi in}{L}x}
\end{equation}\\

Now subbing expression \ref{eq:6} into equation \ref{eq:4}, doing some algebraic manipulation and rearranging the terms we get.

\begin{equation}
\label{eq:6}
\centering \displaystyle \dfrac{\partial \hat{\psi}_n}{\partial t}=-i2(\frac{\pi n}{L})^2\hat{\psi}_n
\end{equation}\\

This is again a first order differential equation which has the following solution.

\begin{equation}
\label{eq:7}
\centering \displaystyle \hat{\psi}_n(t)=\hat{\psi}_n(t_0)e^{-i2(\frac{\pi n}{L})^2}t
\end{equation}\\

\subsection{Finite Difference}

The \textit{Finite Difference} method works by approximating the derivatives in the expression with finite differences. In our \textit{PDE} we have $\frac{\partial \psi}{\partial t}$ and $\frac{\partial^2 \psi}{\partial x^2}$ that need to be approximated via finite differences. The way $\frac{\partial \psi}{\partial t}$ is approximated determines what type of \textit{Finite Difference} scheme is used which has various implications with regards to \textit{accuracy}, \textit{stability} and \textit{implementation}.

\begin{enumerate}
  \item 1) The \textit{Forward Difference} is an \textit{explicit} scheme which means that the solution at each point at the latest time level can be expressed through the solutions of the previous time levels. Although this simplifies the \textit{implementation} the scheme suffers from stability issues. In fact for the \textit{PDE} \ref{eq:1} it can be shown that the \textit{Forward Difference} has an exponential growth in error using the\textit{Von Neumann stability analysis}.
  \item The \textit{Backwards Difference} is an \textit{implicit} scheme which means that a system of equations has to be solved in order to compute the solution at the next time level which makes the implementation non-trivial. However, this method has the superior property that it does not suffer from stability issues.
  \item The \textit{Central Difference} method has the advantage over both the \textit{Forward Difference} and \textit{Backwards Difference} in regards to the accuracy as the error is of $O(\tau ^2)$ compared to $O(\tau)$ of the other methods. This means that the total error of the \textit{PDE} \ref{eq:1} is of $O(\tau ^2 + h^2)$ where $\tau$ is the time step and $h$ is the space step. This method suffers from the same shortcoming as the \textit{Forward Difference} method: Stability issues.
\end{enumerate}

The \textit{Central Difference} was opted to solve the \textit{PDE} \ref{eq:1} as it has good accuracy and although it suffers from stability issues it is stable for certain parameters for $\tau$ and $h$ shown in section 2.5. This method was chosen over the others schemes as the \textit{Forward Difference} is provably unstable and the \textit{Backwards Difference} was not investigated due to the nontrivial implementation and time overhead of solving a system of equations at every time step.

\subsection{Algorithm}

The \textit{Central Difference} depends on the last two time solutions and hence the first time solution was approximated via the \textit{Forward Difference} method. The \textit{Finite Difference} was henceforth implemented as follows.

\begin{enumerate}
	\item Define the initial solution at $t=0$: $ \psi (x, t_0) = f(x)$
	\item Approximate the solution at $t=1$ using the \textit{Forward Difference}: $\psi(x, t_1) = FD( \psi (x, t_0))$
	\item Now the \textit{Central Difference} can be deployed to approximate the solution at $t=i$: $\psi(x, t_i)$ = $CD( \psi (x, t_{i-1}), \psi (x, t_{i-2}))$
\end{enumerate}

where $FD$ and $CD$ are the \textit{Finite Difference} and \textit{Central Difference} schemes respectively which are defined as:

\begin{equation}
\label{FD}
\centering \displaystyle \psi_{j,k+1} = \tau i (0.5 \psi_{xx} + A_{j,k} \psi_{j,k}) + \psi_{j,k}
\end{equation}

\begin{equation}
\label{CD}
\centering \displaystyle \psi_{j,k+1} = 2 \tau i (0.5 \psi_{xx} + A_{j,k} \psi_{j,k}) + \psi_{j,k-1}
\end{equation}

where

\begin{equation}
\label{psixx}
\centering \displaystyle \psi_{xx} = \frac{\psi_{j-1,k} - 2 \psi_{j,k} + \psi_{j+1,k}}{h^2}
\end{equation}

\begin{equation}
\label{A}
\centering \displaystyle A_{j,k} = \frac{| \psi_{j,k} |^2}{1 + S | \psi_{j,k} |^2}
\end{equation}

\subsection{Stability analysis}

In this section the stability of using the \textit{Central Difference} method for \textit{PDE} \ref{eq:1} is analysed to be able to make good choices for the parameters $\tau$ and $h$. The \textit{Von Neumann stability analysis} \cite{analysis} is used to make sense of the stability.\\

Let $\psi_{j,k}=\alpha ^ k e ^ {i \beta j}$ where $\psi_{j,k}$ is the approximation at time $k$ at point $j$. The right hand side of the equation is an arbitrary \textit{Fourier Mode}. If the coefficient $|\alpha|>1$ this implies the solution has an exponential growing error in time and hence it is required that $|\alpha| \leq 1$. \\

We assume the nonlinear term of \textit{PDE} \ref{eq:1} to be negligible and perform the analysis on the linearised version.

\begin{equation}
\label{linPDE}
\centering \displaystyle i\dfrac{\partial \psi}{\partial t} + \dfrac{1}{2} \dfrac{\partial^2 \psi}{\partial x^2} = 0
\end{equation}\\

which has the \textit{Central Difference} formula

\begin{equation}
\label{linPDE}
\centering \displaystyle \psi_{j,k+1} = \frac{\tau i}{h^2} (\psi_{j-1,k} - 2 \psi_{j,k} + \psi_{j + 1,k}) + \psi_{j,k-1}
\end{equation}\\

Substituting $\psi_{j,k}=\alpha ^ k e ^ {i \beta j}$ into this expression we get

\begin{gather*}
	\alpha ^ {k+1}e^{i \beta j} = \frac{\tau i}{h^2}(\alpha ^ k e ^ {i \beta (j - 1)} - 2 \alpha ^ k e ^ {i \beta j} + \alpha ^ k e ^ {i \beta (j + 1)}) + \alpha ^ {k - 1} e^{i \beta j} \\
	\alpha = \frac{\tau i}{h^2} (e ^ {- i \beta} - 2 + e ^ { i \beta}) + \frac{1}{\alpha}\\
	\alpha ^ 2 - \frac{\tau i}{h^2}(2cos\beta - 2) \alpha - 1 = 0\\
	\alpha ^ 2 + \frac{4 \tau i}{h^2}(sin ^ 2 (\frac{\beta}{2})) \alpha - 1 = 0\\
	\alpha = -\frac{2 \tau i}{h^2}sin^2(\frac{\beta}{2}) \pm \sqrt{-\frac{4 \tau^2}{h^4}sin^4(\frac{\beta}{2})+1}
\end{gather*}

The discriminant can either be positive or negative. It will be sufficient to consider the discriminant to be positive as

\begin{equation}
\label{a2=1}
\centering \displaystyle |\alpha|^2 = (-\frac{2 \tau i}{h^2}sin^2(\frac{\beta}{2})))^2 -\frac{4 \tau^2}{h^4}sin^4(\frac{\beta}{2})+1 = 1
\end{equation}

For the discriminant to be positive we require

\begin{equation}
\label{a2=1}
\centering \displaystyle \frac{4 \tau^2}{h^4}sin^4(\frac{\beta}{2})<1
\end{equation}

As we want this inequality to be satisfied for all $\beta$ and $sin^4(\frac{\beta}{2})$ is bounded between $0$ and $1$ we get the inequality.

\begin{equation}
\label{a2=1}
\centering \displaystyle \frac{4 \tau^2}{h^4}<1
\end{equation}

and hence we require

\begin{equation}
\label{a2=1}
\centering \displaystyle \tau < \frac{h^2}{2}
\end{equation}

\subsection{Experimental Setup}

To make a comparison between both methods we require the following quantity to be conserved

\begin{equation}
\label{eq:8}
\centering \displaystyle N = \int_{-\infty}^{\infty} |\psi|^2dx
\end{equation}\\

for a given accuracy $\epsilon$ for one soliton for a short time. \\

In light of this the following parameters were established: The amount of mesh points $N$ were fixed to be $512$ with a simulation time $T=1$ for both methods. The spacial length $L$ was chosen to be $64$ for the \textit{Split Step} method and $30$ for the \textit{Finite Difference} method. The time step $\tau$ for the \textit{Split Step} method was chosen to be $0.01$ (which also satisfied the stability requirement described in section $2.5$). Simulating the solution for $8$ time steps yielded $N=0.00097$. For the \textit{Finite Difference} method the $\tau$ was chosen to be $0.001$ running the simulation for $8$ steps yielded $N=0.00039$. These provided parameters satisfy an $\epsilon=10^{-3}$. $N$ was computed using the composite trapezoidal rule.\\

All implementation was completed using \textit{Python 3} using \textit{Numpy} for the matrix operations, Fourier Transform and Inverse Fourier Transform. \textit{Matplotlib} was used to generate \textit{2D} and \textit{3D} plots. A shared code base was implemented in a file called \textit{helper.py} which contained methods to create initial solutions, plotting functionality and to compute the quantity $N$. Two additional files were generated which implemented the \textit{Split Step} and \textit{Finite Difference} method respectively. \textit{Jupyter Notebook} was used to test both methods.
All the source code can be found in the \textit{Appendix Section}.\\

In the next section the results of multiple simulations using varying values for $S$ and velocities $v_1$ and $v_2$ are portrayed for both methods.

\section{Results}

\subsection{Results with a small negative $S$ of $-0.1$}

\begin{figure}[t]
    \centering
    \includegraphics[width=9cm,height=9cm,keepaspectratio]{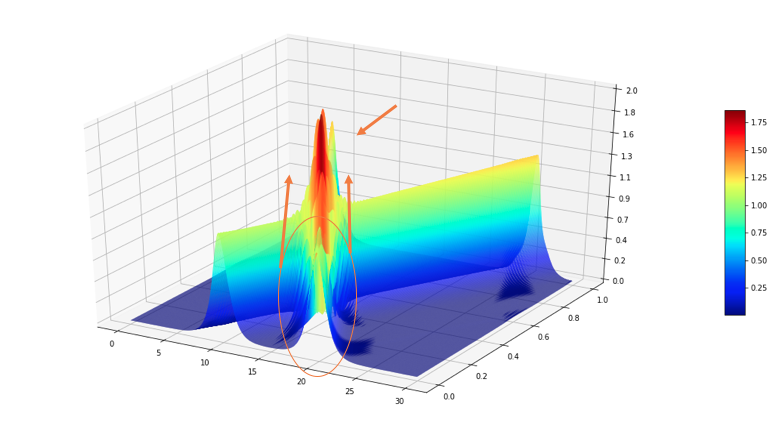}
    \caption{Finite Difference: $S=-0.1$ with $v_1=20$ and $v_2=-20$ (The right hand soliton is a bit tricky to view hence: The oval shows the initial position of the soliton and the arrows attached to the oval point the direction in which it is travelling. The top arrow shows where the soliton is after the collision; Behind the big red spike (superposition))}
    \label{fig1}
\end{figure}

\begin{figure}[t]
    \centering
    \includegraphics[width=9cm,height=9cm,keepaspectratio]{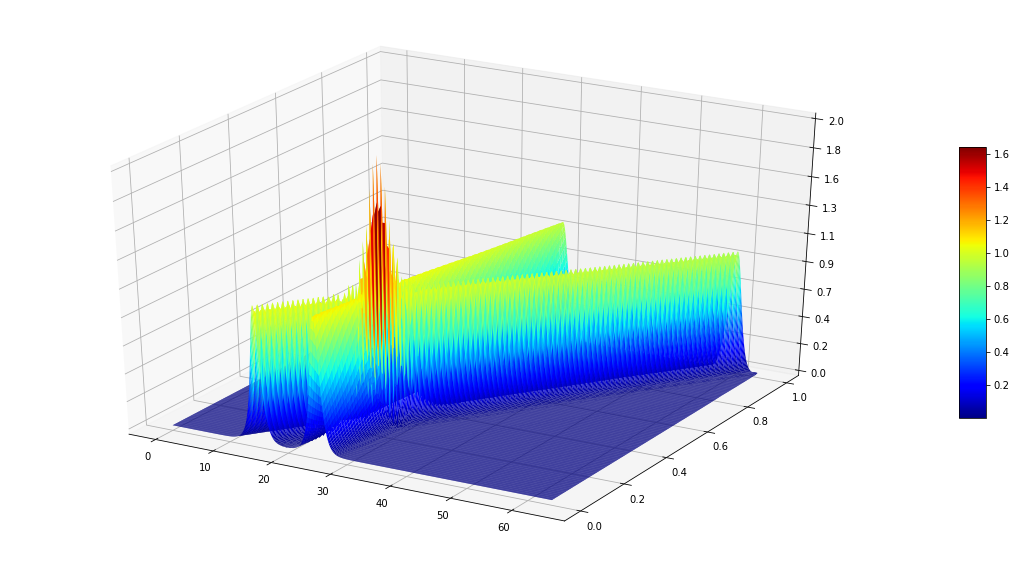}
    \caption{Split Step: $S=-0.1$ with $v_1=20$ and $v_2=-20$}
    \label{fig2}
\end{figure}

\begin{figure}[t]
    \centering
    \includegraphics[width=9cm,height=9cm,keepaspectratio]{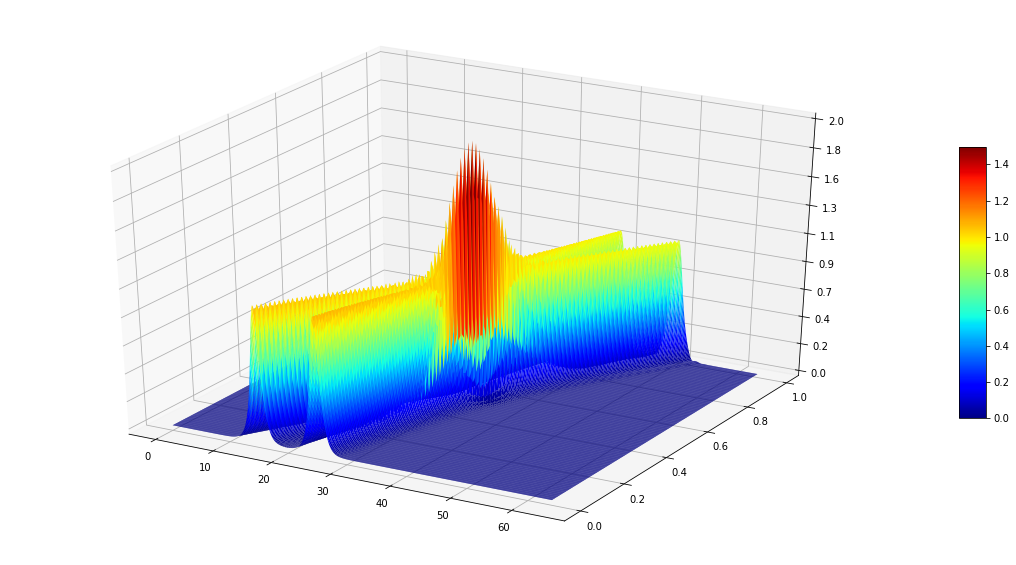}
    \caption{Split Step: $S=-0.1$ with $v_1=10$ and $v_2=-10$}
    \label{fig3}
\end{figure}

\begin{figure}[t]
    \centering
    \includegraphics[width=9cm,height=9cm,keepaspectratio]{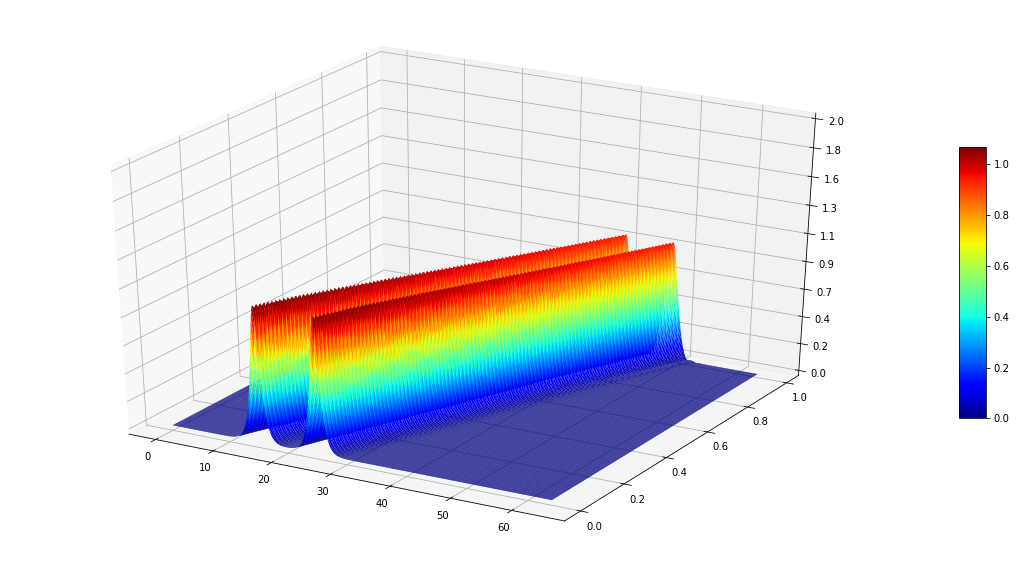}
    \caption{Split Step: $S=-0.1$ with $v_1=1$ and $v_2=-1$}
    \label{fig4}
\end{figure}

For these experiments $S$ was set to a small negative value of $-0.1$. Figure \ref{fig1} is the simulation of the \textit{Finite Difference} method with both solitons initiated with a velocity of $20$ in colliding directions. The same is portrayed in figure \ref{fig2} however using the \textit{Split Step} method. The first observation to be made is that both methods produce colliding solitons that superimpose on each other (\textit{the red spike}) and then decompose back into their original states. A difference between the plots however is that the right soliton of the first figure seems to move left while the right soliton in the next figure seems to rather move in a straight path. \\
Other velocities were also tested using the \textit{Split Step method} (Further experiments using the \textit{Finite Difference} methods were omitted due to time constraints). In figure \ref{fig3} both solitons were initiated with velocities of $10$ in colliding directions. This change seems to be apparent from the visualisation as the solitons collide at a later time step during the simulation. Another experiment was carried out with the solitons initiated with a velocity of $1$ in colliding directions which is depicted in figure \ref{fig4}. In this instance it can be observed that the solitons seems to be travelling parallel to each other as the velocities are not set large enough.

\subsection{Results with a small positive $S$ of $0.4$}

\begin{figure}[t]
    \centering
    \includegraphics[width=9cm,height=9cm,keepaspectratio]{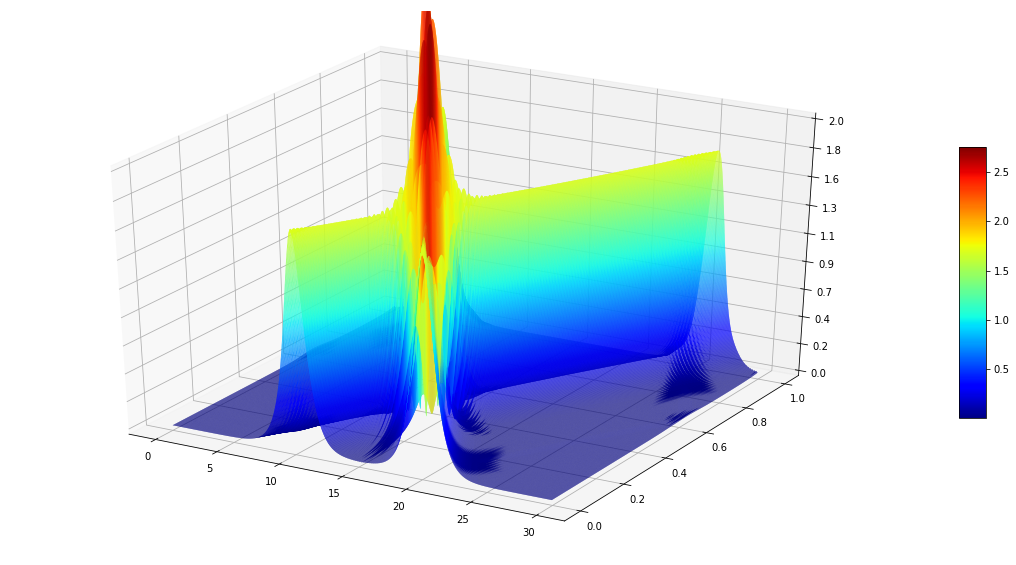}
    \caption{Finite Difference: $S=0.4$ with $v_1=20$ and $v_2=-20$. The right hand solution is hidden behind the big red spike after the the collision. View Figure \ref{fig1} to get a better idea.}
    \label{fig5}
\end{figure}

\begin{figure}[t]
    \centering
    \includegraphics[width=9cm,height=9cm,keepaspectratio]{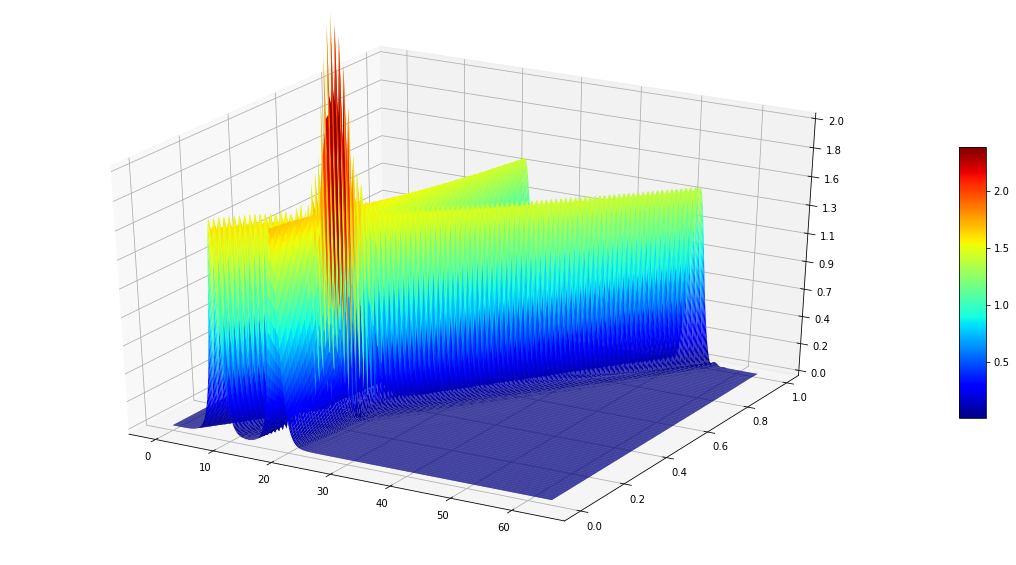}
    \caption{Split Step: $S=0.4$ with $v_1=20$ and $v_2=-20$}
    \label{fig6}
\end{figure}

Two experiments were run, one for the \textit{Finite Difference} and the other for the \textit{Split Step}, however this time with value of $S=0.4$. The velocities were set to $20$ in colliding directions. More changes in velocities were not examined due to time constraints. Figure \ref{fig5} portrays the result using the \textit{Finite Difference}. It can be observed that change in $S$ has caused the solitons to raise in height and breadth. Figure \ref{fig6} shows the same result however using the \textit{Split Step} method from which the same observation can be made.

\subsection{Results with a large negative $S$ of $-10$}

\begin{figure}[t]
    \centering
    \includegraphics[width=9cm,height=9cm,keepaspectratio]{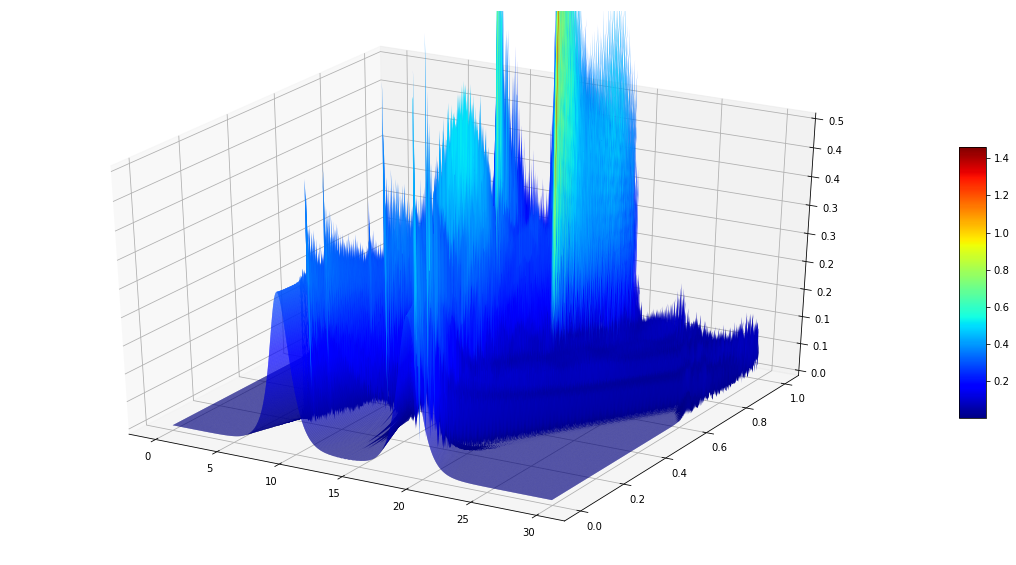}
    \caption{Finite Difference: $S=-10$ with $v_1=10$ and $v_2=-10$}
    \label{fig7}
\end{figure}

\begin{figure}[t]
    \centering
    \includegraphics[width=9cm,height=9cm,keepaspectratio]{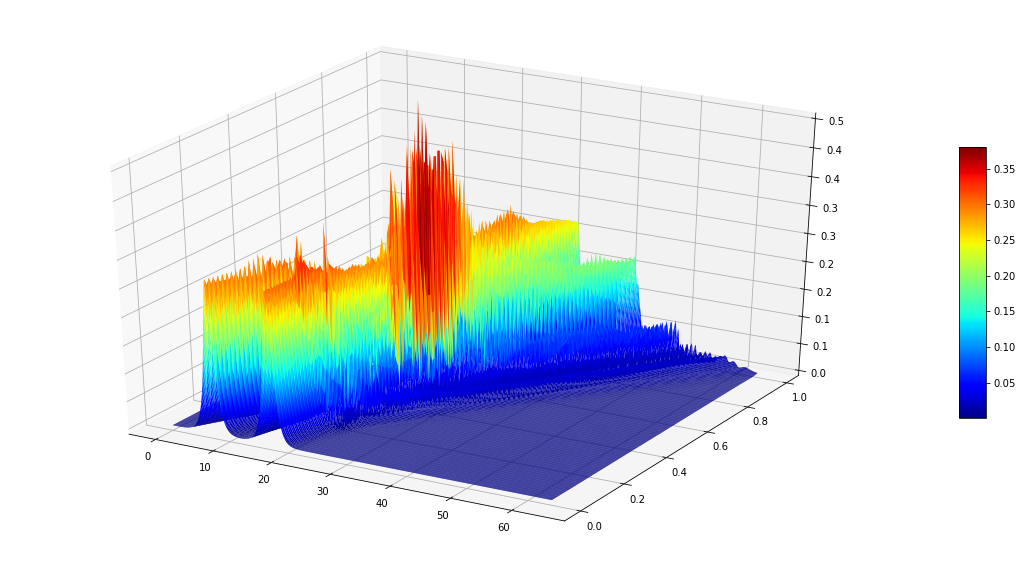}
    \caption{Split Step: $S=-10$ with $v_1=10$ and $v_2=-10$}
    \label{fig8}
\end{figure}

Previous configurations were kept however simulations were now run using a value of $-10$ for $S$ and velocities of $10$ and $-10$ for the two solitons. A few observations can be made: Both initial solutions have a smaller soliton height; Both simulations vary largely from each other. Figure \ref{fig7} results in a \textit{rough} terrain with little structure where else figure \ref{fig8} results in a somewhat more structured output where the solution superposition can be seen, however afterwards small artefact waves can be observed next to the main solitons which seem to dampen after the collision. 

\subsection{Results with a large positive $S$ of $2$}

\begin{figure}[t]
    \centering
    \includegraphics[width=9cm,height=9cm,keepaspectratio]{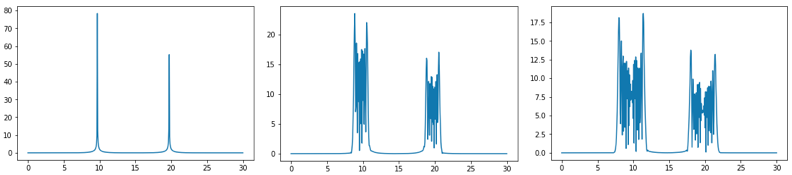}
    \caption{Finite Difference: $S=2$ with $v_1=10$ and $v_2=-10$}
    \label{fig9}
\end{figure}

\begin{figure}[t]
    \centering
    \includegraphics[width=9cm,height=9cm,keepaspectratio]{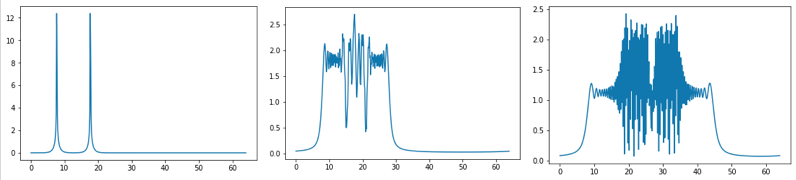}
    \caption{Split Step: $S=2$ with $v_1=10$ and $v_2=-10$}
    \label{fig10}
\end{figure}

Figure \ref{fig9} and figure \ref{fig10} represent simulations run with a value of $S=2$ using the \textit{Finite Difference} and \textit{Split Step} method. The plots provided are \textit{2D} as full simulations were not necessary as from these solutions throughout time it can be observed that the solutions become \textit{noisy} and differ from each other. Increasing $S$ resulted in the the solitons increasing in height and decreasing in width.

\subsection{Run Times}
The time step of the \textit{Split Step} method ran for $5$ seconds and the \textit{Finite Difference} for $2.5$ seconds. However the \textit{Split Step} had a total of $100$ time steps compared to the $1000$ times steps of the \textit{Finite Difference}, hence the total running time using the \textit{Split Step} was $500$ seconds and the running time of the \textit{Finite Difference} was $2500$ seconds which is considerably longer.

\section{Discussion}

In regards to the numeric methods both produced similar plots for values of $S$ close to $0$. A small positive value for $S$ raised the soliton heights where else a small negative value for $S$ reduced the soliton heights. Both schemes produced different results for larger values of $S$ (negative and positive values). This increase in the absolute value of $S$ \textit{increases} the nonlinearity of the \textit{PDE} as it is a coefficient to the absolute value of $\psi$ which has a larger effect on the dynamics of the solutions hence the possible discrepancies between the solutions produced by the two schemes.\\
In regards to the numeric methods the \textit{Split Step} was easier to implement and did not suffer from any instability or careful consideration to choose values $\tau$ and $h$ (in fact further experiments showed that $\tau$ could be increased even more without affecting the accuracy of the outcome). The \textit{Finite Difference} on the other hand was harder to implement and to debug. In addition, stability analysis had to be performed to ensure that the method would actually converge. Finding parameters $\tau$ and $h$ was more difficult compared to the former method. Preliminary experiments demonstrated that even if the stability was satisfied, if $N$ was too small the solutions did not produce desired results. Thus $N$ had to be set large enough yet this forced the step size $\tau$ to be chosen to be considerably smaller than the one chosen for the \textit{Split Step} method. In addition the spacial interval had to be reduced to satisfy the stability condition. \\
The \textit{Split Step} is superior to the \textit{Finite Difference} in this regard as the spacial interval and time step can be relatively large with no stability issues. The \textit{Finite Difference} method on the other hand required careful probing of the these parameters to avoid stability issues. Lastly the \textit{Finite Difference} method requires more computational time than the \textit{Split Step} method due to the small time step.

\section{Conclusion}
In this report the Schrödinger Equation with a particular nonlinearity was investigated using the \textit{Split Step} and \textit{Finite Difference} method. In practice it was found that the \textit{Split Step} method does not suffer from stability issues (like the \textit{Finite Difference}), is faster and can solve the equation on a larger spacial interval using a larger time step. Thus it is advisable to solve problems like these using a spectral method like the \textit{Split Step} over a \textit{Finite Difference method}. Investigation of the equation itself also showed that the parameter $S$ has a large effect on the dynamics of the solutions with smaller values producing solitons that collide with each other, superimpose and then decompose back into their original states whereas larger absolute values destroy the dynamics of the solutions.

\bibliographystyle{named}
\bibliography{ijcai17}

\begin{thebibliography}{}

\bibitem[\protect\citeauthoryear{Delfour \bgroup \em et al.\egroup
  }{1981}]{finitedifference}
M.~Delfour, M.~Fortin, and G.~Payr.
\newblock Finite-difference solutions of a non-linear schrödinger equation.
\newblock In {\em Journal of computational physics}, pages 277--288, 1981.

\bibitem[\protect\citeauthoryear{Kato}{1989}]{nls}
T.~Kato.
\newblock {\em Nonlinear Schrodinger Equations}.
\newblock Springer, Berlin, Germany, 1989.

\bibitem[\protect\citeauthoryear{Keller and Isaacson}{1994}]{analysis}
Herbert Keller and Eugene Isaacson.
\newblock {\em Analysis of numerical methods}.
\newblock Courier Corporation, 1994.

\bibitem[\protect\citeauthoryear{Serkin and Hasegawa}{2000}]{soliton}
V.N. Serkin and A.~Hasegawa.
\newblock Novel soliton solutions of the nonlinear schrödinger equation model.
\newblock In {\em Physical Review Letters}, page 4502, 2000.

\bibitem[\protect\citeauthoryear{Weideman and Herbst}{1986}]{splitstep}
J.A.C. Weideman and B.M. Herbst.
\newblock Split-step methods for the solution of the nonlinear schrödinger
  equation.
\newblock In {\em SIAM Journal on Numerical Analysis}, pages 485--507, 1986.

\bibitem[\protect\citeauthoryear{Zemlyanaya and Alexeeva}{2011}]{nlsoptics}
E.V. Zemlyanaya and N.V. Alexeeva.
\newblock Numerical study of time-periodic solitons in the damped-driven nls.
\newblock In {\em International Journal of Numerical Analysis and Modeling},
  pages 248--261, 2011.

\end{thebibliography}

\section{Appendix}
\subsection{Shared code base: helper.py}
\begin{lstlisting}
import numpy as np
import matplotlib.pyplot as plt
from sympy import *
from mpl_toolkits.mplot3d import Axes3D
import matplotlib.pyplot as plt
from matplotlib import cm
from matplotlib.ticker import LinearLocator, FormatStrFormatter


def initOneSoliton(off, v, L, N, S):
    psi = np.zeros(N, dtype=np.complex_)
    t = 0
    h = L/N
    a = np.sqrt(2)
    B = 3/2 - 2 * S
    for i in range(N):
        x = i * h - off
        f = (2 * a * exp(a * x)) / (1 + B * exp(2 * a * x))
        psi[i] = f * exp(I * t + I * v * x)
    return psi

def initTwoSoliton(x1off, v1, x2off, v2, L, N, S):
    psi = np.zeros(N, dtype=np.complex_)
    t = 0
    h = L/N
    a = np.sqrt(2)
    B = 3/2 - 2 * S
    for i in range(N):
        x1 = i * h - x1off
        f1 = (2 * a * exp(a * x1)) / (1 + B * exp(2 * a * x1))
        x2 = i * h - x2off
        f2 = (2 * a * exp(a * x2)) / (1 + B * exp(2 * a * x2))
        psi[i] = f1 * exp(I * t + I * v1 * x1) + f2 * exp(I * t + I * v2 * x2)
    return psi

def computeN(psi1, psi2, L, N):
    N1 = np.trapz(abs(psi1), dx=L/N)
    N2 = np.trapz(abs(psi2), dx=L/N)
    return np.abs(N1 - N2) 

def plot2D(y, L, N):
    plt.plot(np.linspace(0, L, num=N), abs(y))
    plt.show()

def plot3D(psiEv, L, N, T, tau):
    fig = plt.figure(figsize=(20,10))
    ax = fig.gca(projection='3d')

    # Make data.
    X = np.arange(0, L, L/N)
    Y = np.arange(0, T, tau)
    X, Y = np.meshgrid(X, Y)
    
    # Plot the surface.
    surf = ax.plot_surface(X, Y, psiEv, cmap=cm.jet, linewidth=10, antialiased=True, rstride=1, cstride=1)

    # Customize the z axis.
    ax.set_zlim(0, 2)
    ax.zaxis.set_major_locator(LinearLocator(10))
    ax.zaxis.set_major_formatter(FormatStrFormatter('%.1f'))
    #ax.view_init(30, 190)

    # Add a color bar which maps values to colors.
    fig.colorbar(surf, shrink=0.5, aspect=10)

    plt.show()  
\end{lstlisting}
\subsection{SplitStep.py}
\begin{lstlisting}
import numpy as np
from sympy import *
import helper as h

T = 1
tau = 0.01
L = 64
N = 512
S = -0.1 #3/4 - 0.1

def splitstep(psi):
    # Nonlinear Part
    for i in range(N):
        c = (abs(psi[i]) * abs(psi[i]))
        coef = (I * c) / (1 + S * c)
        psi[i] = psi[i] * exp(coef * tau)
    # Fourier Transform
    c = np.fft.fftshift(np.fft.fft(psi))
    # Move in Fourier Space
    #n = np.linspace(-N/2, N/2+1, num=N)
    for i in range(N):
        e = (-2 * i * i * pi * pi) / (L * L); 
        c[i] = exp(tau * I * e) * c[i]; 
    # Convert back to physical space
    psi = np.fft.ifft(np.fft.fftshift(c)); 
    return psi

# Run the Simulation
psi = h.initTwoSoliton(8, 20, 18, -20, L, N, S)
psiEv = np.zeros(shape=(int(T/tau), int(N)))
for i in range(int(T/tau)):
    psiEv[i] = abs(psi)
    h.plot2D(psi, L, N) # Plot 2D graph at every time step
    psi = splitstep(psi)
h.plot3D(psiEv, L, N, T, tau) # Plot 3D graph

def computeN():
    psi1 = h.initOneSoliton(8, 10, L, N, S)
    psi2 = psi1
    for i in range(8):
        h.plot2D(psi2, L, N)
        psi2 = splitstep(psi2)
    print(h.computeN(psi1, psi2, L, N)) 
\end{lstlisting}
\subsection{FiniteDifference.py}
\begin{lstlisting}
	import numpy as np
from sympy import *
import helper as h

T = 1
tau = 0.001
L = 40
N = 512
S = -0.1 #3/4 - 0.1

def FD(phi1):
    newPhi = np.zeros(N, dtype=np.complex_)
    h = L/N
    
    for i in range(0, N):
        if(i == 0):
            phixx = (phi1[N - 1] - 2 * phi1[i] + phi1[i + 1]) / (h * h)
        elif(i == N-1):
            phixx = (phi1[i - 1] - 2 * phi1[i] + phi1[0]) / (h * h)
        else:
            phixx = (phi1[i - 1] - 2 * phi1[i] + phi1[i + 1]) / (h * h)
        psiSquared = np.abs(phi1[i]) * np.abs(phi1[i])
        A = psiSquared / (1 + S * psiSquared)
        newPhi[i] = I * tau * ((1/2) * phixx + A * phi1[i]) + phi1[i]
        
    return newPhi

def CD(phi1, phi2):
    newPhi = np.zeros(N, dtype=np.complex_)
    h = L/N
    for i in range(0, N):
        if(i == 0):
            phixx = (phi1[N - 1] - 2 * phi1[i] + phi1[i + 1]) / (h * h)
        elif(i == N-1):
            phixx = (phi1[i - 1] - 2 * phi1[i] + phi1[0]) / (h * h)
        else:
            phixx = (phi1[i - 1] - 2 * phi1[i] + phi1[i + 1]) / (h * h)
        psiSquared = np.abs(phi1[i]) * np.abs(phi1[i])
        A = psiSquared / (1 + S * psiSquared)
        newPhi[i] = 2 * I * tau * ((1/2) * phixx + A * phi1[i]) + phi2[i]
    return newPhi

# Run the Simulation
firstPsi = h.initTwoSoliton(10, 20, 20, -20, L, N, S)
secondPsi = FD(firstPsi)
psiEv = np.zeros(shape=(int(T/tau), int(N)))
for i in range(int(T/tau)):
    psiEv[i] = abs(firstPsi)
    h.plot2D(firstPsi, L, N) # Plot 2D graph at every time step
    temp = firstPsi
    firstPsi = secondPsi
    secondPsi = CD(secondPsi, temp)
h.plot3D(psiEv, L, N, T, tau) # Plot 3D graph

def computeN():
    firstPsi = h.initOneSoliton(8, 20, L, N, S)
    secondPsi = FD(firstPsi)
    psi1 = firstPsi
    for i in range(8):
        temp = firstPsi
        firstPsi = secondPsi
        secondPsi = CD(secondPsi, temp)
    print(h.computeN(psi1, secondPsi, L, N)) 
\end{lstlisting}
\end{document}